\documentclass[12pt]{article}
\usepackage{amsfonts}

\newtheorem{thm}{Theorem}[section]
\newtheorem{lemma}[thm]{Lemma}
\newtheorem{cor}[thm]{Corollary}

\newtheorem{prop}[thm]{Proposition}

\newenvironment{remark}{\par\medskip\noindent{\bf Remark.\ }}{\par\smallskip}
\newcommand{\proof
}{\par\medskip\noindent {\bf Proof.\ \ }}

\newcommand{\be}{\begin{equation}}
\newcommand{\ee}{\end{equation}}
\newcommand{\openbox}{\leavevmode
  \hbox to8pt{\hfil\vrule\vbox to6pt{\hrule width6pt\vfil\hrule}\vrule}}

\newcommand{\qed}{\hbox to5pt{ } \hfill \openbox\bigskip\medskip}

\newcommand{\N}{\mathbb N}
\newcommand{\Z}{\mathbb Z}
\newcommand{\Q}{\mathbb Q}

\title{The $f$--vector of the clique complex of chordal graphs
and Betti numbers of edge ideals of uniform hypergraphs}
\author{G\'abor Heged\"{u}s
\\{\normalsize Johann Radon Institute for Computational and Applied Mathematics}
}

\begin{document}


\maketitle

\begin{abstract}
We describe the Betti numbers of the 
edge ideals $I(G)$ of uniform hypergraphs $G$ such that $I(G)$ 
has linear graded free resolution. 

We give an algebraic equation system and some inequalities for the components of 
the $f$--vector
of the clique complex of an arbitrary chordal graph.

Finally we present an explicit formula for the multiplicity of the Stanley-Reisner ring
 of the edge ideals of any chordal graph.
\end{abstract}
\medskip

\section{Introduction}
\noindent

Let $X$ be a finite set and $E:=\{E_1,\ldots ,E_n\}$ a finite collection of non empty 
subsets of $X$. The pair $H=(X,E)$ is called a {\em hypergraph}. The elements of $X$ are 
called the {\em vertices}  and the elements of $E$ are called the {\em edges} of the 
hypergraph. 

We say that a hypergraph $H$ is  {\em $d$-uniform}, if $|E_i|=d$ for every edge 
$E_i\in E$.
 
Let $\Q$ denote the rational field.
Let $R$ be the graded ring $\Q[x_1,\ldots ,x_n]$. 
The vector space $R_s=\Q[x_1,\ldots ,x_n]_s$ consists of 
the homogeneous polynomials of total degree $s$, together with $0$.

We may think of an edge $E_i$ of a hypergraph as a squarefree monomial
 $x^{E_i}:=\prod_{j\in E_i} x_j$ in $R$.

We can associate an ideal $I(H)\subseteq R$ to a hypergraph $H$. 
The {\em edge ideal} $I(H)$ is the ideal $\langle x^{E_i}:~ E_i\in E\rangle$, 
which is generated by the edges of $H$.

The edge ideal was first introduced by R. Villareal in \cite{V}.
Later edge ideals have been studied very widely, see for instance 
\cite{E, E2, F, F2, Fr, HT, HT2, MRV, V, W, Z}.

In \cite{Fr} R. Fr\"oberg characterized the graphs $G$ such that $G$ 
has a linear free resolution. He proved:

\begin{thm} \label{Froberg_theorem}
Let $G$ be  a simple graph on $n$ vertices. Then $R/I(G)$ has linear 
free resolution precisely when $\bar{G}$, the complementary graph of $G$ 
is chordal. 
\end{thm}

In \cite{E2} E. Emtander generalized Theorem \ref{Froberg_theorem} for 
generalized chordal hypergraphs. He proved that the Stanley--Reisner ring of
the incidence complex $\Delta(H)$ corresponding to $H$, 
where $H$ is a generalized chordal hypergraph, has a
 linear free resolution.
In \cite{W} R. Woodroofe extended
 the definition of chordality from graphs to clutters.

In this article we prove explicit formulas for the Betti numbers of the edge
 ideals of $m$-uniform hypergraphs $H$ 
such that $R/I(H)$ has linear free resolution.  

Let $\Delta$ be a simplicial complex. A facet $F$ is called a {\em leaf}, if 
either $F$ is the only facet of $\Delta$, or there exists an other facet
$G$, $G\neq F$ such that $H\cap F\subset G\cap F$ for each facet $H$ with 
$H\neq F$. 
A facet $G$ with this property is called a {\em branch} of $F$. 

Zheng (see \cite{Z}) calls the simplicial complex $\Delta$ a {\em quasi--tree} if
there exists a labeling $F_1,\ldots ,F_m$ of the facets such that 
for all $i$ the facet $F_i$ is a leaf of the subcomplex
 $\langle F_1,\ldots ,F_i\rangle$. We call such a labeling 
a {\em leaf order}. 

A graph is called {\em chordal} if each cycle of length $>3$ has a chord.
 
We recall here for the famous  Dirac's Theorem (see \cite{D}).
\begin{thm} \label{Dirac} (Dirac)
A finite graph $G$ on $[n]$ is a chordal graph iff $G$ is the $1$--skeleton 
of a quasi--tree
\end{thm}

Let $G$ be a finite graph on $[n]$. A {\em clique} of $G$ is a subset 
$F$ of $[n]$ such that $\{i,j\}\in E(G)$ for all $i,j\in F$ with $i\neq j$.

We write $\Gamma(G)$ for the simplicial complex on $[n]$
whose faces are the cliques of $G$.

In our article we give an algebraic equation system 
for the components of the $f$--vector
of the clique complex of an arbitrary chordal graph.

\begin{thm} \label{main_gr} 
Let $G$ be an arbitrary chordal graph. Let $\Gamma:=\Gamma(G)$ be the clique
 complex of
$G$ and $f(\Gamma):=(f_{-1}(\Gamma),\ldots ,f_{d-1}(\Gamma))$ be the 
$f$-vector of the complex $\Gamma$. Here $d=\mbox{dim}(\Gamma)$. Then
\begin{equation} \label{chordal}
-\sum_{i=1}^{p+1} (-1)^i i{f_0 \choose i+1}+ \sum_{j=1}^{p+1} (-1)^{j+p}f_j 
{f_0-j-2\choose p-j+1} =1
\end{equation}
and 
$$
\sum_{k=1}^{p+1} (-1)^k f_k \left( \sum_{i=k-1}^p (-1)^i (2+i)^j {f_0-k-1 \choose i-k+1} \right)+
$$
\begin{equation} \label{chordal2}
+\sum_{i=0}^p (-1)^i (2+i)^j (i+1){f_0\choose i+2}=0,
\end{equation}
for each $j=1,\ldots,n-d-1$, where $p:=pdim(R/I(\overline{G}))$ and $\overline{G}$ is the complement of the graph $G$.
\end{thm}
\begin{remark}
In this Theorem the number of equations 
depends on the dimension of the complex $\Gamma$.
We know from the Auslander--Buchsbaum Theorem that
$n-d\leq p$. If $p=n-d$, then the module $M=R/I(\overline{G})$ is
Cohen--Macaulay and we know that the complement of the chordal graph $G$
is a $d$--tree (see \cite{V2} Theorem 6.7.7, \cite{Fr}). Consequently we know explicitly
the $f$--vector of the clique complices of $d$--trees.
\end{remark}
\begin{thm} \label{main_gr2}
 Let $G$ be an arbitrary chordal graph. Let $\Gamma:=\Gamma(G)$ be the clique
 complex of
$G$ and $f(\Gamma):=(f_{-1}(\Gamma),\ldots ,f_{d-1}(\Gamma))$ be the 
$f$-vector of the complex $\Gamma$. Here $d=\mbox{dim}(\Gamma)$. Then
\begin{equation} \label{chordal4}
\sum_{j=1}^{i+1} (-1)^j f_{j} {f_0-(j+1)\choose i-j+1} +(i+1){f_0 \choose i+2}\geq {p\choose i}
\end{equation}
for each $0\leq i\leq p$, where $p:=pdim(R/I(\overline{G}))$ and $\overline{G}$ is the complement of the graph $G$.
\end{thm}
 
In Section 2 we collected  some basic results about simplicial complices,
 free resolutions, Hilbert fuctions and Hilbert series.
 We present our main results in Section 3.
We prove our main results in Section \ref{proof}.

%%%%%%%%%%%%%%%%%%%%%%%%%%%%%%%%%%%%%%%%%%%%

\section{Preliminaries}

\subsection{Simplicial complices and Stranley--Reisner rings}

We say that $\Delta\subseteq 2^{[n]}$ is a {\em simplicial complex}
 on the vertex set $[n]=\{1,2,\ldots ,n\}$, if 
 $\Delta$ is a set of subsets of $[n]$ such that 
$\Delta$ is a down--set, that is, $G\in\Delta$ and $F\subseteq G$ 
implies that $F\in \Delta$, and $\{i\}\in \Delta$ for all $i$.

The elements of $\Delta$ are called {\em faces} 
and the {\em dimension} of a face is one less than its cardinality. An $r$-face is an abbreviation for an $r$-dimensional face.
The dimension of $\Delta$ is the dimension of a maximal face.
We use the notation $\mbox{dim}(\Delta)$ for the dimension 
of $\Delta$.

If $\mbox{dim}(\Delta)=d-1$, then 
the $(d+1)$--tuple $(f_{-1}(\Delta),\ldots ,f_{d-1}(\Delta))$ 
is called the {\em $f$-vector} of $\Delta$, where
 $f_i(\Delta)$ denotes the number of $i$--dimensional faces 
of $\Delta$.

Let $\Delta$ be an arbitrary simplicial complex
on $[n]$. The {\em Stanley--Reisner ring} 
$R/I_{\Delta}$ of $\Delta$
is the quotient of the ring $R$ by the
 {\em Stanley--Reisner ideal}
$$
I_{\Delta}:=\langle x^F:~ F \notin \Delta \rangle, 
$$ 
generated by the non--faces of $\Delta$.

Let $H=([n], E(H))$ be a simple hypergraph and consider its 
edge ideal $I(H)\subseteq R$. It is easy to verify that $R/I(H)$ is precisely 
the Stanley--Reisner ring of the simplicial complex 
$$
\Delta(H):=\{F\subseteq [n]:~ E\not\subseteq F,\mbox{ for all }E\in E(H)\}.
$$  

This complex is called the {\em independence complex} of $H$. 
By definition the edges of $H$ are precisely the minimal non--faces 
of $\Delta(H)$.

Consider the {\em complementary hypergraph} $\overline{H}$ of a $d$-uniform 
hypergraph. This is defined as the hypergraph $(V(H),E(\bar{H}))$ with
the edge set
$$
E(\overline{H}):=\{F\subseteq X:~ |F|=d,\ F\notin E(H)\}.
$$

Then the edges of $\overline{H}$ are precisely the $(d-1)$-dimensional 
faces of the independence complex $\Delta(H)$. 

Specially, let $H=([n], E(H))$ be a simple graph and consider its 
edge ideal $I(H)\subseteq R$. Then
$$
\Delta(H):=\{F\subseteq [n]:~  F\mbox{ is an independent set in }H\}.
$$  
is the {\em independence complex} of $H$. 
Clearly the edges of $H$ are precisely the minimal non--faces 
of $\Delta(H)$.

Similarly we can define the {\em clique complex} of $H$: 
$$
\Gamma(H):=\{F\subseteq [n]:~ F\mbox{ is a clique in }H\}.
$$ 

\subsection{Free resolutions}

Recall that for every finitely generated graded module $M$ over $R$ 
we can associate to $M$ 
a {\em minimal graded free resolution} 
$$ 
0\longrightarrow \bigoplus_{i=1}^{\beta_p} R(-d_{p,i}) \longrightarrow 
\bigoplus_{i=1}^{\beta_{p-1}} R(-d_{p-1,i})\longrightarrow \ldots \longrightarrow 
\bigoplus_{i=1}^{\beta_{0}} R(-d_{0,i}) \longrightarrow M \longrightarrow
0, 
$$
where $p\leq n$ and $R(-j)$ is the free $R$-module obtained 
by shifting the degrees of $R$ by $j$. 

Here the natural number ${\beta}_{k}$ is the $k$'th {\em total
Betti number} of $M$ and $p$ is the projective dimension of $M$. 

The module $M$ has a {\em pure resolution} if there are 
constants $d_0<\ldots < d_g$
such that 
$$
d_{0,i}=d_0,\ldots ,d_{g,i}=d_g 
$$
for all $i$. If in addition
$$
d_i=d_0+i,
$$
for all $1\leq i\leq p$, then 
we call the minimal free resolution  to be {\em $d_0$--linear}. 

In \cite{R} Theorem 2.7 the following bound for the Betti numbers was proved.
\begin{thm}
Let $M$ be an $R$--module having a pure resolution
of type $(d_0,\ldots, d_p)$ and Betti numbers $\beta_0, \ldots , \beta_p$, where 
$p$ is the projective dimension of $M$.
Then 
\begin{equation} \label{bound_Betti}
\beta_i\geq {p\choose i}
\end{equation}
for each $0\leq i\leq p$.
\end{thm}

\subsection{Hilbert function}

Finally let us recall some basic facts about Hilbert functions and Hilbert series.

Let $M=\bigoplus_{i\geq 0} M_i$ be a finitely generated 
nonnegatively graded module over the polynomial ring $R$. Define 
 the {\em Hilbert function} $h_M : {\Z} \to {\Z}$  by $h_M(i):=\mbox{dim}_{\Q} M_i$.

If we know the $f$-vector of the simplicial complex $\Delta$, 
then we can compute easily the Hilbert function 
$h_{{\Q}[\Delta]}(t)$ of the Stanley--Reisner ring $M:={\Q}[\Delta]$.

\begin{lemma} \label{Stan} (Stanley, see Theorem 5.1.7 in \cite{BH})
The Hilbert function of the Stanley--Reisner ring 
$\Q[\Delta]$ of a $(d-1)$--dimensional simplicial complex $\Delta$ is
\begin{equation} \label{Hilbertfg}
 h_{{\Q}[\Delta]}(t)=\sum_{j=0}^{d-1} f_j(\Delta) {t-1\choose j}.
\end{equation}
\end{lemma} 

In the proof of our main results we use the following Proposition.

\begin{prop} \label{Hilbert} (\cite[Chapter 6, Proposition 4.7]
{CLO2}) Let $M$ be 
a graded $R$-module with the graded free resolution   
\begin{equation}
0\longrightarrow F_n \longrightarrow \ldots \longrightarrow
 F_1 \longrightarrow M \longrightarrow 0.
\end{equation}
If each $F_j$ is the twisted free graded module  $F_j=\bigoplus_{k=1}^{{\beta}_{j,k}} R(d_{j,k})$, then 
\begin{equation}
h_M(t)=\sum_{j=1}^n (-1)^j \sum_{k=1}^{{\beta}_{j,k}} {n+d_{j,k}+t \choose n}.
\end{equation}
\end{prop}

Let $\Delta$ be a simplicial complex such that the Stanley-Reisner 
ring $R/I_{\Delta}$ has a linear free resolution. 
It is known that the generators of $I_{\Delta}$ all have 
the same degree.

It follows that $R/I_{\Delta}$ is a hypergraph algebra $R/I(H)$
for some $k$-uniform hypergraph $H$. 

\subsection{Hilbert--Serre Theorem}

Let $M=\bigoplus_{i\geq 0} M_i$ be a finitely generated 
nonnegatively graded module over the polynomial ring $R$.
We call the formal power series 
$$
H_M(z):= \sum_{i=0}^{\infty} h_M(i)z^i
$$
the {\em Hilbert--series} of the module $M$.

The Theorem of Hilbert--Serre states that there exists 
a (unique) polynomial $P_M(z)\in \Q[z]$, the so-called {\em Hilbert polynomial}
of $M$, such that $h_M(i)=P_M(i)$ for each $i>>0$. Moreover, $P_M$ has degree
$\mbox{dim }M-1$ and $(\mbox{dim }M-1)!$ times the leading coefficient 
of $P_M$ is the {\em multiplicity} of $M$, denoted by $e(M)$. 

Thus, there exist integers $m_0,\ldots, m_{d-1}$ 
such that $h_M(z)=m_0\cdot{z\choose d-1}+m_1\cdot{z\choose d-2}+\ldots + m_{d-1}$,
where ${z\choose r}=\frac{1}{r!}z(z-1)\ldots (z-r+1)$ and $d:=\mbox{dim} M$. 
Clearly $m_0=e(M)$.

We can summarize the Hilbert-Serre theorem as follows: 
\begin{thm} (Hilbert--Serre) \label{Hilbert_Serre}
Let $M$ be a finitely generated nonnegatively graded 
$R$--module of dimension $d$, then the following stetements hold:\\
(a) There exists a (unique) polynomial $P(z)\in \Z[z]$ such that 
the Hilbert--series $H_M(z)$ of $M$ may be written as 
$$
H_M(z)=\frac{P(z)}{(1-z)^d}
$$
(b) $d$ is the least integer for which $(1-z)^dH_M(z)$ is a polynomial.
\end{thm}

%%%%%%%%%%%%%%%%%%%%%%%%%%%%%%%%%%%%%%%%%

\section{The computation of the Betti--vector from the 
$f$-vector}
\subsection{Our main result}

In our main result we describe explicitly the Betti numbers of the 
edge ideals $I(G)$ of uniform hypergraphs $G$ such that $I(G)$ 
has linear free resolution. 

\begin{thm} \label{main6}

Let $G\subseteq {[n] \choose m}$ be an $m$--uniform hypergraph. Suppose 
that the edge ideal $I(G)$ has an $m$-linear free resolution
\begin{equation}
{{\cal F}}_G: 0\longrightarrow R(-m-g)^{\beta_g} \longrightarrow \ldots \longrightarrow 
\end{equation}

\begin{equation} \label{free}
\longrightarrow R(-m-1)^{\beta_1} \longrightarrow R(-m)^{\beta_0} \longrightarrow I(G) \longrightarrow 0.         
\end{equation}

If $\Delta:=\Delta(G)$ is the independence complex of
$G$ and $f(\Delta):=(f_{-1}(\Delta),\ldots ,f_{d-1}(\Delta))$ is the 
$f$-vector of the complex $\Delta$, then

\begin{equation} \label{main3}
\beta_i(G)= \sum_{j=1}^{i+1} (-1)^j f_{j+m-2}(\Delta) {f_0(\Delta)-(j+1)\choose i-j+1} +{i+m-1\choose m-1}{f_0(\Delta) \choose i+m}
\end{equation}
for each $0\leq i\leq g$.
\end{thm}
\begin{remark}
J. Herzog and M. K\"uhl proved similar formulas for the Betti number in \cite{HK}. 
 Theorem 1. Here we did not assume that the ideal $I(G)$ with linear resolution is 
Cohen--Macaulay.   
\end{remark}
\proof
Let $M:=R/I(G)$ denote the quotient module of the edge ideal $I(G)$.
Clearly $R/I(G)$ is the Stanley--Reisner ring of the incidence complex 
$\Delta(G)$.

First we compute the Hilbert function $h_M(t)$ of the 
quotient module $M$ from the graded free resolution of $I(G)$.

From Proposition \ref{Hilbert} 
we conclude that the Hilbert function $h_M(t)$ of $M$ is
\begin{equation}
h_M(t)={t+n\choose n}+\sum_{i=0}^g (-1)^{i+1} \beta_i(G) {t+n-m-i \choose n}.
\end{equation}
From the Vandermonde identities (see e.g. \cite{GKP}, 169--170)
$$
{t+n \choose n}=\sum_{j=0}^{n} {n\choose j}{t\choose j}
$$
and 
$$
{t+n-m-i \choose n}=\sum_{j=0}^{n} {t \choose j}{n-m-i \choose n-j}
$$
for each $i\geq 0$, we infer that
$$
h_M(t)=\sum_{j=0}^n {n\choose j}{t\choose j}+\sum_{i=0}^g (-1)^{i+1} \beta_i(G) \Big( \sum_{j=0}^n {t\choose j} {n-m-i \choose n-j}
\Big)=
$$
$$
=\sum_{j=0}^g {n\choose j}{t\choose j}+
\sum_{j=0}^n {t\choose j} \Big( \sum_{i=0}^g (-1)^{i+1} {n-m-i \choose n-j}\beta_i(G) \Big)
$$
\begin{equation} \label{Hfg1}
=\sum_{j=0}^n {t\choose j} \Big( {n\choose j} +
 \sum_{i=0}^g (-1)^{i+1} {n-m-i \choose n-j}\beta_i(G) \Big)
\end{equation}
On the other hand we can apply Lemma \ref{Stan} for the simplicial complex
 $\Delta$. We get
\begin{equation} \label{Hfg2}
 h_M(t)=\sum_{j=0}^n f_{j-1}(\Delta){t\choose j}.
\end{equation}

But the polynomials $\{{t\choose j}:~ j\in \N \}$ constitute a basis of the 
polynomial ring $\Q[t]$.

Hence equations (\ref{Hfg1}) and (\ref{Hfg2})
imply that 
\begin{equation} \label{Hfg3}
f_{j-1}(\Delta)={n\choose j} + \sum_{i=0}^{j-m} (-1)^{i+1} {n-m-i \choose n-j}\beta_i(G)
\end{equation}
for each $0\leq j\leq n$.

%We can easily solve this linear equation system (\ref{Hfg3}) for $\beta_i(G)$
%to finish the proof of Theorem \ref{main}.

Now we can prove equation (\ref{main3}) by induction. 

It is clear that 
$$
\beta_0(G) +f_{m-1}(\Delta)={n\choose m}.
$$

Hence we settled the case $i=0$.

Suppose that equation (\ref{main3}) is true for each $0\leq i\leq j-m-1$. 
Now we prove equation (\ref{main3}) for $j-m$. 

It follows from equation (\ref{Hfg3}) that
\begin{equation} \label{Hfg4}
(-1)^{j-m}\beta_{j-m}(G)=\sum_{i=0}^{j-m-1} (-1)^{i+1} {n-m-i \choose n-j}\beta_i(G)+{n\choose j}-f_{j-1}(\Delta).
\end{equation}
Hence substituting equation (\ref{main3}) for $\beta_i$, where $0\leq i\leq j-m-1$, 
and rearranging the terms yields to equation (\ref{main3}) for $j-m$.
\qed

In the proof of Theorem \ref{main_gr} we need for the following Corollary.
\begin{cor} \label{main2}

Let $G\subseteq {[n] \choose 2}$ be an $2$--uniform hypergraph. Suppose 
that the edge ideal $I(G)$ has an $2$-linear free resolution
\begin{equation} \label{free3}
{{\cal F}}_G: 0\longrightarrow S(-2-g)^{\beta_g} \longrightarrow \ldots \longrightarrow S(-3)^{\beta_1} \longrightarrow S(-2)^{\beta_0} \longrightarrow I(G) \longrightarrow 0.         
\end{equation}

If $\Delta:=\Delta(G)$ is the independence complex of
$G$ and $f(\Delta):=(f_{-1}(\Delta),\ldots ,f_{d-1}(\Delta))$ is the 
$f$-vector of the complex $\Delta$, then

\begin{equation} \label{main4}
\beta_i(G)= \sum_{j=1}^{i+1} (-1)^j f_{j}(\Delta) {f_0(\Delta)-(j+1)\choose i-j+1} +(i+1){f_0(\Delta) \choose i+2}
\end{equation}
for each $0\leq i\leq g$.
\end{cor}

\subsection{Examples}

We give here two applications of Corollary \ref{main2}.

S. Jacques proved in \cite{J} that
the total $i$'th Betti numbers of the complete graph $K_n$ with $n$ vertices are
$$
\beta_i=(i+1){n\choose i+2}
$$ 
for each $0\leq i\leq n-2$. This is clear from Corollary \ref{main2}, because then
$\bar{G}=([n],\emptyset)$ and the graph $\bar{G}$ is chordal.

Now consider the computation of the total Betti numbers of the 
complete bipartite graphs $K_{n,m}$. Clearly $\overline{K_{n,m}}$ is a chordal graph,
hence it follows from Theorem \ref{Froberg_theorem} that
the edge ideal $I$ has a linear free resolution.

Define the ideal 
$$
I:=I(K_{n,m})=\langle x_iy_j:~ 1\leq i\leq n, 1\leq j\leq m\rangle.
$$

It is easy to verify that the incidence complex $\Delta(K_{n,m})$ is 
the disjoint union of two simplices, one of dimension $n-1$, the other of dimension 
$m-1$.

Hence we get that
$$
f_i(\Delta(K_{n,m}))={n\choose i+1}+{m\choose i+1}
$$
for each $i\geq 0$.

Finally it follows from \cite[Corollary 5.2.5]{J} and Corollary \ref{main2} that
$$
\beta_i(K_{n,m})=\sum_{j+l=i+2,\ j,l\geq 1} {n\choose j}{m\choose l}=
$$
$$
=\sum_{j=1}^{i+1} (-1)^j  \left({n\choose j+1}+{m\choose j+1}\right)
{n+m-j-1\choose i-j+1} +(i+1){n+m \choose i+2}.
$$

%%%%%%%%%%%%%%%%%%%%%%%%%%%%%%%%%%

\section{The proof of our main result}
\subsection{A generalization of Herzog--K\"uhl Theorem}

We need for the following easy Lemma:
\begin{lemma} \label{divis}
Let $K(z)=\sum_{i=0}^p c_iz^{d_i} \in \Q[z]$ be an arbitrary polynomial over $\Q$. Then
$K$ is divisible by $(1-z)^m$ iff $K^{(j)}(1)=0$ for each $j=0,\ldots ,m-1$.    
\end{lemma}

We can prove Theorem \ref{main_gr} with the following generalization of the famous Herzog--K\"uhl
Theorem (Theorem 1 in \cite{HK}). We can prove this Theorem using the same method 
as in \cite{HK}, but for the reader's convenience we include here the proof.
\begin{thm} \label{main_thm}
Let $M$ be an $R$--module having a pure resolution
of type $(d_0,\ldots, d_p)$ and Betti numbers $\beta_0, \ldots , \beta_p$, where 
$p$ is the projective dimension of $M$. Let $d$ denote the dimension of the module 
$M$. Suppose that $d+1\leq n$. 
Then 
\begin{equation} \label{Betti1}
\sum_{i=0}^p (-1)^i \beta_i=0 
\end{equation}
and
\begin{equation} \label{Betti2}
\sum_{i=0}^p (-1)^i \beta_i d_i(d_i-1)\cdot \ldots \cdot (d_i-j+1)=0
\end{equation}
for each $j=1,\ldots , n-d-1$.
\end{thm}
\proof

Since the Hilbert--series is additive on short exact sequences,
and since 
$$
H_R(z)=\frac{1}{(1-z)^n},
$$
and consequently
$$
H_{R(-d)}(z)=\frac{z^d}{(1-z)^n},
$$
the pure resolution
$$ 
0\longrightarrow \bigoplus_{k=1}^{\beta_p} R(-d_{p}) \longrightarrow 
\bigoplus_{k=1}^{\beta_{p-1}} R(-d_{p-1})\longrightarrow \ldots \longrightarrow 
\bigoplus_{k=1}^{\beta_{0}} R(-d_{0}) \longrightarrow M \longrightarrow
0, 
$$
yields 
\begin{equation} \label{Hilbert1}
H_M(z)=\sum_{i=0}^p (-1)^i \beta_i \frac{z^{d_i}}{(1-z)^n}, 
\end{equation}
where $p=pdim(M)$.

Write $d:=\mbox{dim} M$, and let $m:=\mbox{codim}(M)=n-d$. 
It follows from the Auslander--Buchbaum formula that $m\leq p$. 
We infer from the 
Theorem of Hilbert--Serre that we can write 
\begin{equation} \label{Hilbert2}
H_M(z)=\frac{P(z)}{(1-z)^d}.
\end{equation}
Comparing the two expressions (\ref{Hilbert1}) and (\ref{Hilbert2})
for $H_M$, we find
\begin{equation} \label{main}
(1-z)^m P(z)=\sum_{i=0}^p (-1)^i \beta_i z^{d_i}
\end{equation}
This formula shows that $(1-z)^m$ divides $\sum_{i=0}^p (-1)^i \beta_i z^{d_i}$ (in the 
ring $\Z[x]$). It follows from Lemma \ref{divis} that 
$(\beta_0,\ldots , \beta_p)$ solves the equation system (\ref{Betti1}), (\ref{Betti2}).  \qed

\subsection{The multiplicity of Stanley-Reisner ideals of chordal graphs}

We can derive easily the following Corollary.
\begin{cor}
Let $M$ be an $R$--module having a pure resolution
of type $(d_0,\ldots, d_p)$ and Betti numbers $\beta_0, \ldots , \beta_p$, where 
$p$ is the projective dimension of $M$. Let $d$ denote the dimension of the module 
$M$. Suppose that $d+1\leq n$. 
Then 
\begin{equation} \label{Betti4}
\sum_{i=0}^p (-1)^i \beta_i d_i^j=0
\end{equation}
for each $j=0,\ldots , n-d-1$.
\end{cor}
\begin{remark}
It follows easily that these equations are linearly independent. 
\end{remark}
\begin{cor} \label{mult}
Let $M$ be an $R$--module having a pure resolution
of type $(d_0,\ldots, d_p)$ and Betti numbers $\beta_0, \ldots , \beta_p$, where 
$p$ is the projective dimension of $M$. Let $d$ denote the dimension of the module 
$M$ and $m:=\mbox{codim}(M)=n-d$. Suppose that $d+1\leq n$. 
Then 
$$
e(M)=(-1)^m\frac{p!}{m!} \sum_{i=0}^p (-1)^i \beta_i{d_i\choose p}.
$$
\end{cor}
\proof
It comes out from the definition that
$$
e(M)=\left((1-z)^d\cdot H_M(z) \right)\mid_{z=1}=P(1).
$$
Hence we infer from equation (\ref{main}) that
$$
e(M)=P(1)=\frac{(-1)^m}{m!}\left((1-z)^mP \right)^{(m)}\mid_{z=1}= 
$$
$$
=\frac{(-1)^m}{m!} \sum_{i=0}^p (-1)^i \beta_i p!{d_i\choose p}=
$$
$$
=(-1)^m\frac{p!}{m!} \sum_{i=0}^p (-1)^i \beta_i{d_i\choose p}.
$$
\qed

Now we can describe easily the multiplicity of the Stanley--Reisner 
ideals of chordal graphs.
\begin{cor}
Let $G$ be an arbitrary
 chordal graph and $H:=\overline{G}$ denote the complement of the graph $G$.
Let $\Gamma:=\Gamma(G)$ be the clique
 complex of
$G$ and $f(\Gamma):=(f_{-1}(\Gamma),\ldots ,f_{d-1}(\Gamma))$ be the 
$f$-vector of the complex $\Gamma$. Let $p$ be the projective dimension of $R/I(H)$. 
 Let $d$ denote the dimension of the module $R/I(H)$ and $m:=\mbox{codim}(R/I(H))=n-d$. 
Then
$$
e(R/I(H))=(-1)^m\frac{p!}{m!} \sum_{i=0}^p (-1)^i \left( \sum_{j=1}^{i+1} (-1)^j f_{j} {f_0-(j+1)\choose i-j+1} +(i+1){f_0 \choose i+2}\right) {i+2\choose p}
$$
\end{cor}
\proof 
It follows from Theorem \ref{Froberg_theorem} that 
the module $M:=R/I(H)$ has a $2$-linear resolution:
\begin{equation} \label{free4}
{{\cal F}}_H: 0\longrightarrow S(-2-p)^{\beta_p} \longrightarrow \ldots \longrightarrow S(-3)^{\beta_1} \longrightarrow S(-2)^{\beta_0} \longrightarrow R \longrightarrow M \longrightarrow 0.         
\end{equation}
where $p$ is the projective dimension of $M$.

If we apply Theorem \ref{mult} for the module $M$, we get that
\begin{equation} \label{mult2}
e(R/I(H))=(-1)^m\frac{p!}{m!} \sum_{i=0}^p (-1)^i \beta_i{i+2\choose p}.
\end{equation}
Now using Theorem \ref{main2} and substituting 
$$
\beta_i(H)= \sum_{j=1}^{i+1} (-1)^j f_{j} {f_0-(j+1)\choose i-j+1} +(i+1){f_0 \choose i+2}
$$
into (\ref{mult2}), we get our result. \qed

\subsection{The proofs} \label{proof}

{\bf Proof of Theorem \ref{main_gr}:}
Let $H:=\overline{G}$ denote the complement of the graph $G$.

Then consider the module $M:=R/I(H)$.
It follows from Theorem \ref{Froberg_theorem} that 
the module $M$ has a $2$-linear resolution:
\begin{equation} \label{free2}
{{\cal F}}_H: 0\longrightarrow S(-2-p)^{\beta_p} \longrightarrow \ldots \longrightarrow S(-3)^{\beta_1} \longrightarrow S(-2)^{\beta_0} \longrightarrow R \longrightarrow M \longrightarrow 0.         
\end{equation}
where $p$ is the projective dimension of $M$.

If we apply Theorem \ref{main2} for the graph $F:=H$, then we get that
\begin{equation} \label{main5}
\beta_i(H)= \sum_{j=1}^{i+1} (-1)^j f_{j}(\Delta) {f_0(\Delta)-(j+1)\choose i-j+1} +(i+1){f_0(\Delta) \choose i+2}
\end{equation}
for each $0\leq i\leq p$.

Now we can apply Theorem \ref{main_thm}. If we substitute the expressions (\ref{main5}) for $\beta_i(H)$ into the equation system
(\ref{Betti1}), (\ref{Betti2}) and rearrange the obtained equations, we get 
our result. 

Namely
$$
\sum_{i=0}^p (-1)^i\beta_i=\sum_{i=0}^p (-1)^i \left(\sum_{j=1}^{i+1} (-1)^j 
f_j{f_0-j-1 \choose i-j+1}+(i+1){f_0 \choose i+2} \right)
$$
$$
=\sum_{i=0}^p (-1)^i (i+1){f_0 \choose i+2}+\sum_{i=0}^p (-1)^i
 \left( \sum_{j=1}^{i+1} (-1)^j f_j{f_0-(j+1)\choose i-j+1} \right)
$$
$$
=\sum_{i=0}^p (-1)^i (i+1){f_0 \choose i+2}+ \sum_{j=1}^{p+1} (-1)^jf_j(\Gamma) 
\left( \sum_{i=j-1}^p (-1)^i{f_0-(j+1)\choose i-j+1} \right)
$$
$$
=-\sum_{i=1}^{p+1} (-1)^i i{f_0 \choose i+1}+ \sum_{j=1}^{p+1} (-1)^{j+p}f_j(\Gamma) 
{f_0(\Gamma)-j-2\choose p-j+1} =-1,
$$
because 
$$
\sum_{i=j-1}^p (-1)^i{f_0-(j+1)\choose i-j+1} = (-1)^p{f_0-j-2 \choose p-j+1}.
$$
Similarly 
$$
\sum_{i=0}^p (-1)^id_i^j \beta_i=\sum_{i=0}^p (-1)^i (2+i)^j
\left( \sum_{j=1}^{i+1} (-1)^j f_j{f_0-j-1 \choose i-j+1}+(i+1){f_0 \choose i+2} \right)
$$
$$
\sum_{k=1}^{p+1} (-1)^k f_k(\Gamma) \left( \sum_{i=k-1}^p (-1)^i (2+i)^j {f_0-k-1 \choose i-k+1} \right)+
$$
\begin{equation} \label{chordal3}
+\sum_{i=0}^p (-1)^i (2+i)^j (i+1){f_0\choose i+2}=0.
\end{equation} 
\qed

{\bf Proof of Theorem \ref{main_gr2}}:
Let $H:=\overline{G}$ denote the complement of the graph $G$. Applying
Theorem \ref{Froberg_theorem} and Theorem \ref{main2} for the graph $F:=H$,
we get again (\ref{main5}). Hence we infer from Theorem \ref{bound_Betti} that
\begin{equation} \label{chordal5}
\sum_{j=1}^{i+1} (-1)^j f_{j} {f_0-(j+1)\choose i-j+1} +(i+1){f_0 \choose i+2}\geq {p\choose i}
\end{equation}
for each $0\leq i\leq p$. \qed

%\begin{equation} \label{chordal3}
%\sum_{i=0}^p (-1)^i (i+1){f_0 \choose i+2}+ \sum_{j=1}^{p+1} (-1)^jf_j(\Gamma) 
%\left( \sum_{i=j-1}^p (-1)^i{f_0-(j+1)\choose i-j+1} \right)=-1
%\end{equation}
%$$
%\sum_{k=1}^{p+1} (-1)^k f_k(\Gamma) \left( \sum_{i=k-1}^p (-1)^i (2+i)^j {f_0-k-1 \choose i-k+1} \right)+
%$$
%\begin{equation} \label{chordal4}
%+\sum_{i=0}^p (-1)^i (2+i)^j (i+1){f_0\choose i+2}=0,
%\end{equation}
{\bf Acknowledgements.}  I am indebted to Josef Schicho and Lajos R\'onyai 
for their useful remarks.


\begin{thebibliography}{MM}
\bibitem{BH} W. Bruns, J. Herzog, {\em Cohen-Macaulay rings},
 Cambridge Studies in Advanced Mathematics, {\bf 39},
 Cambridge University Press, Cambridge, 1993. 

\bibitem{CLO} D. Cox, J. Little, D. O'Shea, {\em Ideals, varieties, and
algorithms}, Springer, 1992. 

\bibitem{CLO2} D. Cox, J. Little, D. O'Shea, {\em Using Algebraic Geometry}, Springer, 2005.

\bibitem{D} Dirac, G. A., On rigid circuit graphs.
{\em Abh. Math. Sem. Univ. Hamburg} {\bf 25} (1961) 71–76. 

\bibitem{E} E. Emtander, Betti numbers of hypergraphs.
  {\em Comm. Algebra}  {\bf 37}  (2009),  no. 5, 1545--1571.

\bibitem{E2} E. Emtander,A class of hypergraphs that generalizes chordal graphs. 
{\em  Math. Scand.}  {\bf 106}  (2010),  no. 1, 50--66. 

\bibitem{F} S. Faridi, Cohen-Macaulay properties of square-free monomial ideals. {\em 
 J. Combin. Theory Ser.} A  {\bf 109}  (2005),  no. 2, 299--329. 

\bibitem{F2} S. Faridi,  The facet ideal of a simplicial complex 
{\em Manuscripta Math.} {\bf 109} (2002), no. 2, 159--174. 

\bibitem{Fr} R. Fr\"oberg, On Stanley-Reisner rings.  Topics in algebra, Part 2 (Warsaw, 1988),  
57--70, Banach Center Publ., 26, Part 2, PWN, Warsaw, 1990.

\bibitem{GKP} R.L. Graham, D.E. Knuth, O. Patashnik, {\em Concrete
Mathematics,} Addison-Wesley, 1989.

\bibitem{HT} H. T. H$\grave{a}$; A. Van Tuyl, Monomial ideals, edge ideals of
 hypergraphs, and their graded Betti numbers.  {\em J. Algebraic Combin.} 
 {\bf 27}  (2008),  no. 2, 215--245.

\bibitem{HT2} H. T. H$\grave{a}$; A. Van Tuyl,
Splittable ideals and the resolutions of monomial ideals. 
 {\em J. Algebra} {\bf 309} (2007),  no. 1, 405--425.

\bibitem{HHZ} J. Herzog,  T. Hibi, X. Zheng,
Dirac's theorem on chordal graphs and Alexander duality. 
{\em European J. Combin.} {\bf 25} (2004), no. 7, 949–960. 

\bibitem{H} G. Heged\"us, Betti numbers of edge ideals of uniform hypergraphs, 
see {\em http://arxiv.org/abs/1009.0394}

\bibitem{HK} J. Herzog,  M. K\"uhl,
On the Betti numbers of finite pure and linear resolutions.
{\em Comm. Algebra} {\bf 12} (1984), no. 13-14, 1627–1646. 

\bibitem{J} S. Jacques, Betti Numbers of Graph Ideals, PhD Thesis, available online 
'http://arxiv.org/abs/math/0410107'

\bibitem{MS} E. Miller; B. Sturmfels, Combinatorial commutative algebra. 
Graduate Texts in Mathematics, 227. Springer-Verlag, New York, 2005.

\bibitem{MRV} S. Morey; E. Reyes; R. H. Villarreal,
 Cohen-Macaulay, shellable and unmixed clutters with a perfect 
matching of König type. {\em  J. Pure Appl. Algebra } {\bf 212}  (2008),  no. 7, 1770--1786

\bibitem{R} T. R\"omer, Bounds for Betti numbers. {\em J. Algebra} {\bf 249} (2002), no. 1, 20–-37.

\bibitem{V} R. H. Villarreal,  Cohen-Macaulay graphs.  {\em Manuscripta Math.} 
 {\bf 66}  (1990),  no. 3, 277--293. 

\bibitem{V2} R. H.Villarreal,  Monomial algebras.
 Monographs and Textbooks in Pure and Applied Mathematics
, {\bf 238}. Marcel Dekker, Inc., New York, 2001

\bibitem{W} R. Woodroofe, Chordal and sequentially Cohen-Macaulay clutters, available online 
'http://arxiv.org/abs/0911.4697'

\bibitem{Z} Zheng, X. Resolutions of facet ideals. 
{\em Comm. Algebra} {\bf 32} (2004), no. 6, 2301–2324.
\end{thebibliography}
\end{document}